\documentclass{ifacconf}

\usepackage{graphicx}      
\usepackage{natbib}        

\usepackage{amsmath}
\usepackage{amssymb}

\begin{document}
\begin{frontmatter}

\title{System-theoretic Analysis of Nonlinear Infinite-dimensional Systems with Generalized Symmetries\thanksref{footnoteinfo}} 

\thanks[footnoteinfo]{This work has been supported by the Austrian Science Fund (FWF) under grant number P 29964-N32.}

\author[First]{Bernd Kolar} 
\author[First]{Markus Sch{\"o}berl} 

\address[First]{Institute of Automatic Control and Control Systems Technology, Johannes Kepler University, Linz, Austria \\ (e-mail: bernd.kolar@jku.at, markus.schoeberl@jku.at)}

\end{frontmatter}

\setlength{\parindent}{15pt}
\setlength{\parskip}{0pt}

\section{Introduction}

Symmetry groups are groups that transform solutions of partial or ordinary differential equations continuously into other solutions. Classical symmetry groups act geometrically on the space of independent and dependent variables of the differential equation, and are generated by vector fields on this space. For this reason, they are also called geometric symmetries or point symmetries.
Generalized symmetry groups, in contrast, are generated by generalized vector fields,
which include derivatives of the dependent variables,
and act directly on the infinite-dimensional space of solutions. Therefore, as already recognized by E. Noether, they significantly extend the applicability of symmetry group methods, see e.g. \cite{Olver:1993}.

A special case of classical symmetry groups are the vertical symmetry groups, which transform only the dependent variables. In e.g. \cite{RiegerSchlacher:2010} and \cite{KolarRamsSchoberl:2018}, it has been shown that these groups can be useful for analyzing the system-theoretic property observability of infinite-dimensional systems with input and output. More precisely, the existence of certain vertical symmetry groups implies the non-observability of the system.
However, the focus on vertical symmetry groups is a severe restriction. Most symmetry groups of nonlinear systems also transform the independent variables, see e.g. \cite{Ibragimov:1994}, and for systems that do not possess vertical symmetries the method is not applicable.
In this contribution, we discuss how the existing approach based on vertical classical symmetries can be extended to non-vertical classical symmetries, as well as to arbitrary generalized symmetries.

\section{Generalized Symmetries of PDEs with Boundary Conditions}

Vertical symmetry groups are generated by vertical vector fields, i.e., vector fields which have no components in direction of the independent variables. For non-vertical symmetry groups, the evolutionary form of the corresponding non-vertical vector fields plays an important role. The reason is that we have to determine the influence of the symmetry group on functions that are defined at fixed points of the spatial domain, like boundary conditions or point outputs. For this purpose, we have to use the Lie derivative along the evolutionary form of the vector field. To illustrate this, consider a system of PDEs with independent variables $z,t$ and dependent variables $x^\alpha$, $\alpha=1,\ldots,q$. A vector field on the space of independent and dependent variables has the form
\begin{equation}
 v=v_z(z,t,x)\partial_z+v_t(z,t,x)\partial_t+v_x^\alpha(z,t,x)\partial_{x^\alpha}\,,
 \label{eq:VectorField}
\end{equation}
and its evolutionary form is the (vertical) generalized vector field
\begin{equation}
 v_Q=(v_x^\alpha(z,t,x)-v_z(z,t,x)x_z^\alpha-v_t(z,t,x)x_t^\alpha)\partial_{x^\alpha}\,,
 \label{eq:EvolutionaryForm}
\end{equation}
which also depends on the first derivatives of $x$ with respect to $z$ and $t$. For vertical vector fields, (\ref{eq:VectorField}) and (\ref{eq:EvolutionaryForm}) coincide.
Now suppose we have a solution $x^\alpha=\gamma^\alpha(z,t)$ of the PDE, and the symmetry group generated by the vector field (\ref{eq:VectorField}) transforms it for every (sufficiently small) value of the group parameter $\varepsilon$ into a new solution $x^\alpha=\gamma_\varepsilon^\alpha(z,t)$. Then, for a smooth function $c(x)$, the change of the composed function $c\circ \gamma_\varepsilon(z,t)$ caused by a variation of the group parameter $\varepsilon$ is given by
\begin{equation}
 \left. \partial_\varepsilon \left( c\circ \gamma_\varepsilon(z,t) \right) \right|_{\varepsilon=0} = L_{v_Q}c \circ j^1(\gamma(z,t))\,.
 \label{eq:LvQ_c}
\end{equation} 
To take account of the components of $v$ in $z$- and $t$-direction, here we need the Lie derivative of $c(x)$ along the evolutionary form $v_Q$ of the vector field, and not along the vector field itself. For functions $c(x,x_z,x_t,\ldots)$ that depend also on derivatives of $x$, the Lie derivative must be taken along a suitable prolongation of $v_Q$.
Both the formula (\ref{eq:EvolutionaryForm}) for the calculation of the evolutionary form as well as the relation (\ref{eq:LvQ_c}) remain valid for generalized vector fields $v$, with functions $v_z$, $v_t$, and $v_x^\alpha$ that depend also on derivatives of $x$ up to some order. Depending on the order of these derivatives, it is of course necessary to replace $j^1(\gamma(z,t))$ in (\ref{eq:LvQ_c}) by a sufficiently high prolongation.

Let the system of PDEs now be given by
\begin{equation}
 \Delta^\nu=0\,,\quad \nu=1,\ldots,l\,,
 \label{eq:Delta}
\end{equation}
where $\Delta^\nu$ are functions on an appropriate jet space over the bundle with coordinates $(z,t,x)\rightarrow (z,t)$. Furthermore, assume that the Jacobian matrix of the functions $\Delta^\nu$ has full rank, which guarantees that the equations (\ref{eq:Delta}) describe a regular submanifold of this jet space. According to \cite{Olver:1993}, a generalized vector field $v$ is a generalized infinitesimal symmetry if the Lie derivatives
\begin{equation}
 L_{j(v)}\Delta^\nu\,,\quad \nu=1,\ldots,l
 \label{eq:Lie_Delta}
\end{equation}
vanish for every smooth solution $x=\gamma(z,t)$ of the PDE. Here $j(v)$ denotes a sufficiently high prolongation of $v$. As shown in \cite{Olver:1993}, Proposition 5.5, it does not matter whether we take the Lie derivative in (\ref{eq:Lie_Delta}) along the prolongation of $v$ or the prolongation of $v_Q$.

In control theory, PDEs are usually accompanied by boundary conditions. Thus, we must ensure additionally that the symmetry group transforms solutions that satisfy the boundary conditions into solutions that also satisfy the boundary conditions. Geometrically, boundary conditions can be written as equations
\begin{equation}
 \Delta_{BC}^\lambda=0\,,\quad \lambda=1,\ldots,m
 \label{eq:DeltaBC}
\end{equation}
on the restriction of appropriate jet spaces to the respective parts of the boundary. Hence, the generalized vector field must meet the additional condition that the Lie derivatives
\begin{equation}
 L_{j(v_Q)}\Delta_{BC}^\lambda\,,\quad \lambda=1,\ldots,m
 \label{eq:Lie_DeltaBC}
\end{equation}
vanish for every smooth solution $x=\gamma(z,t)$ that satisfies the boundary conditions. In contrast to (\ref{eq:Lie_Delta}), here it is essential to use the evolutionary form $v_Q$. Again, $j(v_Q)$ denotes a sufficiently high prolongation. If the (generalized) vector field meets both (\ref{eq:Lie_Delta}) and (\ref{eq:Lie_DeltaBC}), it is a (generalized) infinitesimal symmetry of the system with boundary conditions. It should be noted that this condition is sufficient but not necessary, since we could relax (\ref{eq:Lie_Delta}) by requiring only that the Lie derivatives vanish for every smooth solution that satisfies the boundary conditions. Unfortunately, this criterion is much harder to check.

\section{Observability Analysis}

A pair of initial conditions of a nonlinear system with input and output is said to be indistinguishable, if for every
admissible trajectory of the input, the system generates for both initial conditions the same trajectory of the output. If there exists no pair of indistinguishable initial conditions the system is said to be observable, otherwise it is not observable.
In the following, we consider for simplicity like in \cite{RiegerSchlacher:2010} only the case of a fixed input trajectory, which leads to an autonomous, time-variant system. In this simplified setting, two initial conditions are indistinguishable if they generate the same output trajectory.
The main idea is now just like in \cite{RiegerSchlacher:2010} and \cite{KolarRamsSchoberl:2018}, but instead of vertical classical symmetry groups we make use of generalized symmetry groups: If there exists a generalized symmetry of the system which does not change the trajectory of the output, then it can be used to transform solutions into other solutions with indistinguishable initial conditions. Thus, the system cannot be observable.
In the following, the approach is demonstrated by an example which does not possess a vertical classical symmetry group.

Consider the autonomous system
\begin{equation}
 \partial_t x(z,t)=(x(z,t)+1)\partial_z x(z,t)
 \label{eq:Example_PDE}
\end{equation}
on the spatial domain $\Omega=(0,1)$ with boundary condition
\begin{equation}
 x(1,t)=0
 \label{eq:Example_BC}
\end{equation} and a point output
\begin{equation}
 y(t)=\left. \tfrac{\partial_z x(z,t)}{x(z,t)} \right|_{z=0}\,.
 \label{eq:Example_y}
\end{equation}This system is similar to the elementary nonlinear wave equation discussed in \cite{Olver:1993}, but extended by the boundary condition and the output function.
It can be shown easily that the vector field
\begin{equation}
 v=zx\partial_z+(x+1)x\partial_x
 \label{eq:Example_v}
\end{equation}
with evolutionary form
\begin{equation}
 v_Q=(x+1-zx_z)x\partial_x
 \label{eq:Example_vQ}
\end{equation}
is an infinitesimal symmetry of the PDE (\ref{eq:Example_PDE}) with boundary condition (\ref{eq:Example_BC}). The Lie derivative
\begin{equation}
 L_{j^1(v)}\Delta
\end{equation}
of the function $\Delta=x_t-(x+1)x_z$ representing the PDE (\ref{eq:Example_PDE}) vanishes for all solutions of the PDE, and the Lie derivative
\begin{equation}
 \left. L_{v_Q}\Delta_{BC}\right|_{z=1} = (x+1-x_z)x
\end{equation}
of the function $\Delta_{BC}=x$ representing the boundary condition (\ref{eq:Example_BC}) vanishes for all solutions that satisfy the boundary condition, since at the boundary $z=1$ we have $x=0$. Thus, the vector field is an infinitesimal symmetry of the system (\ref{eq:Example_PDE}) with boundary condition (\ref{eq:Example_BC}). Furthermore, also the Lie derivative
\begin{equation}
 \left. L_{j^1(v_Q)}\tfrac{x_z}{x}\right|_{z=0}=\left. -zx_{zz}\right|_{z=0} = 0
\end{equation}
of the output function (\ref{eq:Example_y}) vanishes. Thus, the group action generated by the vector field (\ref{eq:Example_v}) maps solutions onto solutions with the same output trajectory $y(t)$ (provided that the group action indeed exists). Hence, the initial conditions of these solutions are indistinguishable, and the system is not observable.

The calculation of the group action is more difficult than in the case of vertical classical symmetry groups. Given a solution $x(z,t)=\gamma (z,t)$, the transformed solution $x(z,t)=\gamma_\varepsilon(z,t)$ for some value of the group parameter $\varepsilon$ is determined by the solution of the Cauchy problem
\begin{equation}
 \partial_\varepsilon x(z,t,\varepsilon)=(x(z,t,\varepsilon)+1-z\partial_z x(z,t,\varepsilon))x(z,t,\varepsilon)
 \label{eq:Cauchy_Problem}
\end{equation}
with boundary condition $x(1,t,\varepsilon)=0$ and initial condition $x(z,t,0)=\gamma (z,t)$. The right-hand side of (\ref{eq:Cauchy_Problem}) is the coefficient of the evolutionary vector field (\ref{eq:Example_vQ}), and the boundary condition is due to the fact that the transformed solution must satisfy the boundary condition (\ref{eq:Example_BC}) for all values of the group parameter $\varepsilon$.
For a rigorous proof of the non-observability of the system, we would have to calculate the group action, or at least show that it indeed exists. However, as pointed out in \cite{Olver:1993}, proving the existence and uniqueness of solutions of such Cauchy problems is a very difficult problem, which lies far beyond the scope of this paper. Therefore, we confine ourselves to the verification of the infinitesimal conditions.

\section{Conclusion}

The extension from vertical classical symmetries to non-vertical classical symmetries and generalized symmetries makes the observability analysis suggested in \cite{RiegerSchlacher:2010} and \cite{KolarRamsSchoberl:2018} applicable to a much larger class of systems, since many nonlinear systems do not possess vertical classical symmetries at all.
However, infinitesimal generators of generalized symmetries are harder to find, and the calculation of the group action from the infinitesimal generator is a difficult task.


\begin{thebibliography}{4}
\providecommand{\natexlab}[1]{#1}
\providecommand{\url}[1]{\texttt{#1}}
\providecommand{\urlprefix}{URL }
\expandafter\ifx\csname urlstyle\endcsname\relax
  \providecommand{\doi}[1]{doi:\discretionary{}{}{}#1}\else
  \providecommand{\doi}{doi:\discretionary{}{}{}\begingroup
  \urlstyle{rm}\Url}\fi

\bibitem[{Ibragimov(1994)}]{Ibragimov:1994}
Ibragimov, N. (ed.) (1994).
\newblock \emph{CRC Handbook of Lie Group Analysis of Differential Equations},
  volume~1.
\newblock CRC Press, Boca Raton.

\bibitem[{Kolar et~al.(2018)Kolar, Rams, and
  Sch{\"o}berl}]{KolarRamsSchoberl:2018}
Kolar, B., Rams, H., and Sch{\"o}berl, M. (2018).
\newblock Application of symmetry groups to the observability analysis of
  partial differential equations.
\newblock In \emph{Proceedings 23rd International Symposium on Mathematical
  Theory of Networks and Systems (MTNS)}, 247--254.

\bibitem[{Olver(1993)}]{Olver:1993}
Olver, P. (1993).
\newblock \emph{Applications of Lie Groups to Differential Equations}.
\newblock Springer, New York, 2nd edition.

\bibitem[{Rieger and Schlacher(2010)}]{RiegerSchlacher:2010}
Rieger, K. and Schlacher, K. (2010).
\newblock Accessibility and observability for a class of first-order {PDE}
  systems with boundary control and observation.
\newblock In \emph{Proceedings 19th International Symposium on Mathematical
  Theory of Networks and Systems (MTNS)}.

\end{thebibliography}
\end{document}